\documentclass[english]{article}
\usepackage[T1]{fontenc}
\usepackage[latin9]{inputenc}
\usepackage{amstext}
\usepackage{amssymb}
\usepackage{esint}

\makeatletter

\begin{document}

\title{Reducing a class of two-dimensional integrals to one-dimension with
application to Gaussian Transforms}

\maketitle
\markright{Reducing a class of 2D integrals ti 1D}

\author{Jack C. Straton}
\begin{abstract}
Quantum theory is awash in multidimensional integrals that contain exponentials in the integration variables, their inverses, and inverse polynomials of those variables. The present paper introduces a means to reduce pairs of such integrals to one dimension when the integrand contains powers times an arbitrary function of $xy/(x+y)$ multiplying various combinations of exponentials. In some cases these exponentials arise directly from transition-amplitudes involving products of plane waves, hydrogenic wave functions, Yukawa and/or Coulomb potentials. In other cases these exponentials arise from Gaussian transforms of such functions. 
\end{abstract}

\section{Introduction}

Prudnikov, Brychkov, and Marichev \cite{PBM p. .567 No. 3.1.3.4}
provide a means to reduce a class of two-dimensional integrals involving
exponentials and a general function of the integration variables in
the specific configuration $f\left(\frac{xy}{x+y}\right)$ to a one-dimensional
integral involving that same function of the new integration variable
$f\left(t\right)$,

\begin{eqnarray}
\int_{0}^{\infty}\int_{0}^{\infty}\frac{1}{\sqrt{x+y}} & f & \left(\frac{xy}{x+y}\right)e^{-px-qy}dx\,dy\nonumber \\
 & = & \frac{\sqrt{\pi}}{\sqrt{pq\left(p+q\right)}}\int_{0}^{\infty}e^{-\left(\sqrt{p}+\sqrt{q}\right)^{2}t}f(t)\,dt\quad.\label{eq:PBM p. .567 No. 3.1.3.4}
\end{eqnarray}

As may be seen with comparison to the specific cases, such as entry
No. 3.1.3.5 on the same page, the coefficient is in error. This reduction
should read

\begin{eqnarray}
\int_{0}^{\infty}\int_{0}^{\infty}\frac{1}{\sqrt{x+y}} & f & \left(\frac{xy}{x+y}\right)e^{-px-qy}dx\,dy\nonumber \\
 & = & \frac{\sqrt{\pi}\left(\sqrt{p}+\sqrt{q}\right)}{\sqrt{p}\sqrt{q}}\int_{0}^{\infty}e^{-\left(\sqrt{p}+\sqrt{q}\right)^{2}t}f(t)\,dt\quad.\label{eq:fixed}
\end{eqnarray}

\noindent In the present paper we extend this reduction technique
to a class of integrals that arise in Gaussian transforms of atomic,
molecular, and optical transition amplitudes \cite{Stra89b,Stra90a}
in which the above exponentials on the left-hand side may contain
more complicated forms, $e^{-\frac{a}{x}-\frac{b}{y}-c\,xy/(x+y)-h\,y/(x+y)-j/(x+y)-px-qy}$,
and one may also have positive and negative powers of x and y.

We will begin with a more general form

\vspace{0.5cm}

\begin{eqnarray}
R_{2} & \left(\right. & n\left.,\,m,\,\nu,\,a,\,b,\,c,\,h,\,j,\,p,\,q\right)=\int_{0}^{\infty}\int_{0}^{\infty}\frac{1}{x^{n/2}y^{m/2}\left(x+y\right)^{\nu/2}}f\left(\frac{xy}{x+y}\right)\nonumber \\
 & \times & e^{-\frac{a}{x}-\frac{b}{y}-c\,xy/(x+y)-h\,y/(x+y)-j/(x+y)-px-qy}dx\,dy\quad.\label{eq:I2gen}
\end{eqnarray}

\noindent The reduction is facilitated by a simultaneous change of
variables to

\begin{eqnarray}
s=\frac{x}{\sqrt{x+y}} & , & t=\frac{sy}{\sqrt{x+y}}\\
x=s^{2}+t\equiv\phi\left(s,\,t\right) & , & y=\frac{t\left(s^{2}+t\right)}{s^{2}}\equiv\psi\left(s,\,t\right)\\
\frac{\partial\phi}{\partial s}\frac{\partial\psi}{\partial t}-\frac{\partial\psi}{\partial s}\frac{\partial\phi}{\partial t} & = & \frac{2\left(s^{2}+t\right)^{2}}{s^{3}}\quad,
\end{eqnarray}
 where the last line is the Jacobian determinant of the transformation.
We retain the integration limits over $\left[0,\infty\right]$. Then

\begin{eqnarray}
R_{2} & \left(\right. & n\left.,\,m,\,\nu,\,a,\,b,\,c,\,h,\,j,\,p,\,q\right)\nonumber \\
 & = & 2\int_{0}^{\infty}f(t)\,t^{-m/2}e^{-ct-qt}\,dt\int_{0}^{\infty}s^{m+\nu-3}\left(s^{2}+t\right)^{-\frac{m}{2}-\frac{n}{2}-\nu+2}\label{eq:I2stgen}\\
 & \times & \exp\left(-\frac{a}{s^{2}+t}-\frac{bs^{2}}{t\left(s^{2}+t\right)}-\frac{ht}{s^{2}+t}-\frac{js^{2}}{\left(s^{2}+t\right)^{2}}-p\left(s^{2}+t\right)-\frac{qt^{2}}{s^{2}}\right)ds\nonumber 
\end{eqnarray}

\section{Set of integral reductions for positive powers in the exponentials }

\noindent If we set $a=b=c=h=j=0$, complete the square in the exponential
and change variables in the $s$ integral

\begin{equation}
J_{2}\left(n,\,m,\,\nu,\,p,\,q\right)=\int_{0}^{\infty}s^{m+\nu-3}\left(s^{2}+t\right)^{-\frac{m}{2}-\frac{n}{2}-\nu+2}\exp\left(-ps^{2}-\frac{qt^{2}}{s^{2}}\right)ds\label{eq:J2gens}
\end{equation}
 to

\begin{eqnarray}
u & = & \sqrt{p}s-\frac{\sqrt{q}t}{s},s=\frac{\sqrt{4\sqrt{p}\sqrt{q}t+u^{2}}+u}{2\sqrt{p}}\\
ds & = & \frac{\left(\sqrt{4\sqrt{p}\sqrt{q}t+u^{2}}+u\right)^{2}}{\sqrt{p}\left(\sqrt{4\sqrt{p}\sqrt{q}t+u^{2}}+u\right)^{2}+4p\sqrt{q}t}\text{du}
\end{eqnarray}
 then

\vspace{0.7cm}

\begin{eqnarray}
J_{2} & \left(\right. & n\left.,\,m,\,\nu,\,p,\,q\right)=e^{-2\sqrt{p}\sqrt{q}t}\int_{-\infty}^{\infty}e^{-u^{2}}\label{eq:J2gen}\\
 & \times & \frac{2^{-m-\nu+3}p^{\frac{1}{2}(-m-\nu+3)}\left(\frac{\left(\sqrt{4\sqrt{p}\sqrt{q}t+u^{2}}+u\right)^{2}}{4p}+t\right)^{-\frac{m}{2}-\frac{n}{2}-\nu+2}}{\left(\sqrt{p}\left(\sqrt{4\sqrt{p}\sqrt{q}t+u^{2}}+u\right)^{2}+4p\sqrt{q}t\right)\left(\sqrt{4\sqrt{p}\sqrt{q}t+u^{2}}+u\right)^{1-m-\nu}}\,du\nonumber 
\end{eqnarray}

\subsection{Integrals with $m+\nu=1$}

\noindent If we set $\{n=0,m=0,\nu=1\}$, the second factor in the
denominator of (\ref{eq:J2gen}) goes to unity, allowing us to do
the integral \cite{GR5 p. 385 No. 3.472.3}, and we recover the corrected
integral (\ref{eq:fixed}) from Prudnikov, Brychkov, and Marichev
\cite{PBM p. .567 No. 3.1.3.4}

\begin{eqnarray}
R_{2} & \left(\right. & 0\left.,\,0,\,1,\,0,\,0,\,0,\,0,\,0,\,p,\,q\right)=\int_{0}^{\infty}\int_{0}^{\infty}\frac{1}{\sqrt{x+y}}f\left(\frac{xy}{x+y}\right)e^{-px-qy}dx\,dy\nonumber \\
 & = & \frac{\sqrt{\pi}\left(\sqrt{p}+\sqrt{q}\right)}{\sqrt{p}\sqrt{q}}\int_{0}^{\infty}f(t)\,e^{-\left(\sqrt{p}+\sqrt{q}\right)^{2}t}\,dt\quad.\label{eq:corrected PBM p. .567 No. 3.1.3.4}
\end{eqnarray}

\noindent A second form may be obtained by setting 
\begin{equation}
f\left(\frac{xy}{x+y}\right)\equiv\left(\frac{xy}{x+y}\right)^{-1/2}g\left(\frac{xy}{x+y}\right)\quad\label{eq:f_mp5tog}
\end{equation}
 in the first form, giving

\begin{eqnarray}
R_{2} & \left(\right. & 1\left.,\,1,\,0,\,0,\,0,\,0,\,0,\,0,\,p,\,q\right)=\int_{0}^{\infty}\int_{0}^{\infty}\frac{1}{x^{1/2}y^{1/2}}g\left(\frac{xy}{x+y}\right)e^{-px-qy}dx\,dy\nonumber \\
 & = & \frac{\sqrt{\pi}\left(\sqrt{p}+\sqrt{q}\right)}{\sqrt{p}\sqrt{q}}\int_{0}^{\infty}t^{-1/2}g(t)\,e^{-\left(\sqrt{p}+\sqrt{q}\right)^{2}t}\,dt\quad.\label{eq:mnn=00003D00003D00003D110}
\end{eqnarray}

\noindent Additional applications of (\ref{eq:f_mp5tog}) will give
larger powers of \emph{x} and \emph{y} in the denominator with powers
of $x+y$ in the numerator. Applications of the inverse of (\ref{eq:f_mp5tog})
will give powers of \emph{x} and \emph{y} in the numerator with larger
powers of $x+y$ in the denominator.

One may also do the integral if we create even powers for the last
factor in the numerator of (\ref{eq:J2gen}) and expand, such as setting
$\{n=-1,m=1,\nu=0\}$, which gives dissimilar powers of \emph{x} and
\emph{y} \cite{Mathematica7},

\begin{eqnarray}
R{}_{2} & \left(\right. & -1\left.,\,1,\,0,\,0,\,0,\,0,\,0,\,0,\,p,\,q\right)=\int_{0}^{\infty}\int_{0}^{\infty}\frac{x^{1/2}}{y^{1/2}}f\left(\frac{xy}{x+y}\right)e^{-px-qy}dx\,dy\nonumber \\
 & = & \frac{\sqrt{\pi}}{2p^{3/2}\sqrt{q}}\int_{0}^{\infty}f(t)\,e^{-\left(\sqrt{p}+\sqrt{q}\right)^{2}t}\left(2\sqrt{p}t\left(\sqrt{p}+\sqrt{q}\right)^{2}+\sqrt{q}\right)\,dt\quad.\label{eq:Im110}
\end{eqnarray}

\noindent The transformation (\ref{eq:f_mp5tog}) can be applied to
this and all subsequent integrals.

\subsection{Integrals with $n+m+2\nu=4$}

\noindent We obtain a new set of integral reductions by finding values
of \emph{n, m}, and $\nu$ that set the last factor in the numerator
of (\ref{eq:J2gen}) to unity.

\subsubsection{Integrals with $\nu=0$}

\noindent For $\{n=3,m=1,\nu=0\}$, and using the identity

\begin{equation}
f_{1}\left(\frac{xy}{x+y}\right)=\left(\frac{xy}{x+y}\right)g\left(\frac{xy}{x+y}\right)\quad,\label{eq:f_1tog}
\end{equation}
 we obtain another set of integral reductions for unlike powers of
the original coordinates

\begin{eqnarray}
R_{2} & \left(\right. & 3\left.,\,1,\,0,\,0,\,0,\,0,\,0,\,0,\,p,\,q\right)=R_{2}\left(1,\,-1,\,2,\,0,\,0,\,0,\,0,\,0,\,p,\,q\right)\label{eq:I2w310}\\
 & = & \int_{0}^{\infty}\int_{0}^{\infty}\frac{\sqrt{y}}{\sqrt{x}(x+y)}g\left(\frac{xy}{x+y}\right)e^{-px-qy}dx\,dy\nonumber \\
 & = & \frac{\sqrt{\pi}}{\sqrt{q}}\int_{0}^{\infty}g(t)\,t^{-1/2}e^{-\left(\sqrt{p}+\sqrt{q}\right)^{2}t}\,dt\quad\left[\Re\left(\sqrt{p}\sqrt{q}\right)>0\right]\quad.\nonumber 
\end{eqnarray}
 For $\{n=5,m=-1,\nu=0\}$ we employ (\ref{eq:f_2tog}) to obtain

\begin{eqnarray}
R_{2} & \left(\right. & 5\left.,\,-1,\,0,\,0,\,0,\,0,\,0,\,0,\,p,\,q\right)=R_{2}\left(1,\,-5,\,4,\,0,\,0,\,0,\,0,\,0,\,p,\,q\right)\label{eq:I2w5m10}\\
 & = & \int_{0}^{\infty}\int_{0}^{\infty}\frac{y^{5/2}}{\sqrt{x}(x+y)^{2}}g\left(\frac{xy}{x+y}\right)e^{-px-qy}dx\,dy\nonumber \\
 & = & \frac{\sqrt{\pi}}{2q^{3/2}}\int_{0}^{\infty}g(t)\,t^{-1/2}\left(1+2\sqrt{p}\sqrt{q}t\right)e^{-\left(\sqrt{p}+\sqrt{q}\right)^{2}t}\,dt\quad\left[\Re\left(\sqrt{p}\sqrt{q}\right)>0\right]\quad.\nonumber 
\end{eqnarray}
 One may continue on in like manner.

\subsubsection{Integrals yielding Macdonald functions}

\noindent One can also set $\{n=1,m=1,\nu=1\}$ to set the $\left(s^{2}+t\right)$
term of (\ref{eq:J2gens}) to unity to obtain \cite{GR5 p. 384 No. 3.471.9}

\begin{eqnarray}
R_{2} & \left(\right. & 1\left.,\,1,\,1,\,0,\,0,\,0,\,0,\,0,\,p,\,q\right)\nonumber \\
 & = & \int_{0}^{\infty}\int_{0}^{\infty}\frac{1}{x^{1/2}y^{1/2}\left(x+y\right)^{1/2}}f\left(\frac{xy}{x+y}\right)e^{-px-qy}dx\,dy\nonumber \\
 & = & 2\int_{0}^{\infty}f(t)\,t^{-1/2}e^{-pt-qt}K_{0}\left(2\sqrt{p}\sqrt{q}t\right)\,dt\quad\left[\Re\left(\sqrt{p}\sqrt{q}\right)>0\right].\label{eq:I2w111}
\end{eqnarray}

\noindent We can instead integrate \cite{Mathematica7} the $u$-form
with the last factor in the numerator of (\ref{eq:J2gen}) set to
one with $\{n=1,m=1,\nu=1\}$. Since 
\[
J_{2}\left(1,\,1,\,1,\,p,\,q\right)=J_{2}\left(2,\,2,\,0,\,p,\,q\right)
\]
 we also have,

\begin{eqnarray}
R_{2} & \left(\right. & 2\left.,\,2,\,0,\,0,\,0,\,0,\,0,\,0,\,p,\,q\right)=\int_{0}^{\infty}\int_{0}^{\infty}\frac{1}{xy}f\left(\frac{xy}{x+y}\right)e^{-px-qy}dx\,dy\nonumber \\
 & = & 2\int_{0}^{\infty}f(t)\,t^{-1}e^{-pt-qt}K_{0}\left(2\sqrt{p}\sqrt{q}t\right)\,dt\quad\left[\Re\left(\sqrt{p}\sqrt{q}\right)>0\right],\label{eq:I2w220}
\end{eqnarray}
 but this may also be obtained from (\ref{eq:I2w111}) by using (\ref{eq:f_mp5tog}).

We can also set the last factor in the numerator of (\ref{eq:J2gen})
to unity with $\{n=4,m=0,\nu=0\}$ to obtain a set of integral reductions
that have unlike powers of the original coordinates \cite{GR5 p. 384 No. 3.471.9}
,

\begin{eqnarray}
R_{2} & \left(\right. & 4\left.,\,0,\,0,\,0,\,0,\,0,\,0,\,0,\,p,\,q\right)=\int_{0}^{\infty}\int_{0}^{\infty}\frac{1}{x^{2}}f_{2}\left(\frac{xy}{x+y}\right)e^{-px-qy}dx\,dy\nonumber \\
 & = & 2\frac{\sqrt{p}}{\sqrt{q}}\int_{0}^{\infty}f_{2}(t)\,t^{-1}e^{-pt-qt}K_{1}\left(2\sqrt{p}\sqrt{q}t\right)\,dt\quad\left[\Re\left(\sqrt{p}\sqrt{q}\right)>0\right]\quad.\label{eq:I2w400}
\end{eqnarray}
 It is to be understood that each of the integral reductions in this
paper are valid only for those functions\emph{ f} that are convergent
in the final integral. In the above equation we have made this restriction
more explicit with the notation $f_{2}(t)$ that indicates that the
integrals do not converge for terms in a powers series representation
of $f$ for powers less than 2. But if we set

\begin{equation}
f_{2}\left(\frac{xy}{x+y}\right)=\left(\frac{xy}{x+y}\right)^{2}g\left(\frac{xy}{x+y}\right)\quad.\label{eq:f_2tog}
\end{equation}
 then we can rewrite

\begin{eqnarray}
R_{2} & \left(\right. & 4\left.,\,0,\,0,\,0,\,0,\,0,\,0,\,0,\,p,\,q\right)=R_{2}\left(0,\,-4,\,4,\,0,\,0,\,0,\,0,\,0,\,p,\,q\right)\nonumber \\
 & = & \int_{0}^{\infty}\int_{0}^{\infty}\frac{y^{2}}{(x+y)^{2}}g\left(\frac{xy}{x+y}\right)e^{-px-qy}dx\,dy\label{eq:I2w400-1}\\
 & = & 2\frac{\sqrt{p}}{\sqrt{q}}\int_{0}^{\infty}g(t)\,t^{-1+2}e^{-pt-qt}K_{1}\left(2\sqrt{p}\sqrt{q}t\right)\,dt\quad\left[\Re\left(\sqrt{p}\sqrt{q}\right)>0\right],\nonumber 
\end{eqnarray}
 a form that converges for any function \emph{g }that can be expanded
in a powers series (without the negative powers of a Laurant series).
For {$n=0,m=4,\nu=0\}$ we simply interchange $p\leftrightarrow q$
in the above.

Similarly, for $\{n=5,n=1,\nu=-1\}$ we employ (\ref{eq:f_2tog})
to obtain

\vspace{-0.06cm}

\begin{eqnarray}
R_{2} & \left(\right. & 5\left.,\,1,\,-1,\,0,\,0,\,0,\,0,\,0,\,p,\,q\right)=R_{2}\left(1,\,-5,\,4,\,0,\,0,\,0,\,0,\,0,\,p,\,q\right)\nonumber \\
 & = & \int_{0}^{\infty}\int_{0}^{\infty}\frac{y^{3/2}}{\sqrt{x}(x+y)^{3/2}}g\left(\frac{xy}{x+y}\right)e^{-px-qy}dx\,dy\label{eq:I2w5m1m1}\\
 & = & 2\frac{\sqrt{p}}{\sqrt{q}}\int_{0}^{\infty}g(t)\,t^{1/2}e^{-pt-qt}K_{1}\left(2\sqrt{p}\sqrt{q}t\right)\,dt\quad\left[\Re\left(\sqrt{p}\sqrt{q}\right)>0\right]\quad.\nonumber 
\end{eqnarray}
 One may continue on in this fashion with more extreme powers such
as $\{n=7,n=-1,\nu=-1\}$, by using the identity

\begin{equation}
f_{3}\left(\frac{xy}{x+y}\right)=\left(\frac{xy}{x+y}\right)^{3}g\left(\frac{xy}{x+y}\right)\quad,\label{eq:f_3tog}
\end{equation}
 to obtain

\begin{eqnarray}
R_{2} & \left(\right. & 7\left.,\,-1,\,-1,\,0,\,0,\,0,\,0,\,0,\,p,\,q\right)=R_{2}\left(1,\,-7,\,5,\,0,\,0,\,0,\,0,\,0,\,p,\,q\right)\nonumber \\
 & = & \int_{0}^{\infty}\int_{0}^{\infty}\frac{y^{7/2}}{\sqrt{x}(x+y)^{5/2}}g\left(\frac{xy}{x+y}\right)e^{-px-qy}dx\,dy\label{eq:I2w400-2}\\
 & = & \frac{2p}{q}\int_{0}^{\infty}g(t)\,t^{3/2}e^{-pt-qt}K_{2}\left(2\sqrt{p}\sqrt{q}t\right)\,dt\quad\left[\Re\left(\sqrt{p}\sqrt{q}\right)>0\right]\quad.\nonumber 
\end{eqnarray}

\subsection{Integrals with $m+\nu=2$}

\noindent We obtain a new set of integral reductions by finding values
of \emph{m} and $\nu$ that move the second factor in the denominator
of (\ref{eq:J2gen}) to the numerator with unit power. The ratio of
this term with the other term left in the denominator reduces nicely
\cite{Mathematica7},

\begin{equation}
\frac{\left(\sqrt{4\sqrt{p}\sqrt{q}t+u^{2}}+u\right)^{1}}{\left(\sqrt{p}\left(\sqrt{4\sqrt{p}\sqrt{q}t+u^{2}}+u\right)^{2}+4p\sqrt{q}t\right)}=\frac{1}{\sqrt{4\sqrt{p}\sqrt{q}t+u^{2}}}\label{eq:ratio_red}
\end{equation}
 allowing integrals of nonzero powers of the numerator of (\ref{eq:J2gen})
to be done. We already saw the results for $\{n=1,m=1,\nu=1\}$ in
(\ref{eq:I2w111}). For $\{n=-1,m=1,\nu=1\}$ we obtain \cite{Mathematica7}

\begin{eqnarray}
R_{2} & \left(\right. & -1\left.,\,1,\,1,\,0,\,0,\,0,\,0,\,0,\,p,\,q\right)\nonumber \\
 & = & \int_{0}^{\infty}\int_{0}^{\infty}\frac{x^{1/2}}{y^{1/2}\left(x+y\right)^{1/2}}f\left(\frac{xy}{x+y}\right)e^{-px-qy}dx\,dy\nonumber \\
 & = & \frac{2}{\sqrt{p}}\int_{0}^{\infty}f(t)\,t^{+1/2}e^{-pt-qt}\left(\sqrt{p}K_{0}\left(2\sqrt{p}\sqrt{q}t\right)+\sqrt{q}K_{1}\left(2\sqrt{p}\sqrt{q}t\right)\right)\,dt\nonumber \\
 &  & \quad\left[\Re\left(\sqrt{p}\sqrt{q}\right)>0\right].\label{eq:I2wm111}
\end{eqnarray}
 For $\{n=-3,m=1,\nu=1\}$ \cite{Mathematica7},

\begin{eqnarray}
R_{2} & \left(\right. & -3\left.,\,1,\,1,\,0,\,0,\,0,\,0,\,0,\,p,\,q\right)\nonumber \\
 & = & \int_{0}^{\infty}\int_{0}^{\infty}\frac{x^{3/2}}{y^{1/2}\left(x+y\right)^{1/2}}f\left(\frac{xy}{x+y}\right)e^{-px-qy}dx\,dy\nonumber \\
 & = & \frac{2}{p^{3/2}}\int_{0}^{\infty}f(t)\,t^{+1/2}e^{-pt-qt}\nonumber \\
 & \times & \left(\sqrt{p}t(p+q)K_{0}\left(2\sqrt{p}\sqrt{q}t\right)+\sqrt{q}(2pt+1)K_{1}\left(2\sqrt{p}\sqrt{q}t\right)\right)\,dt\nonumber \\
 &  & \quad\left[\Re\left(\sqrt{p}\sqrt{q}\right)>0\right]\nonumber \\
 & = & \frac{\partial}{\partial p}R_{2}\left(-1,\,1,\,1,\,0,\,0,\,0,\,0,\,0,\,p,\,q\right)\quad.\label{eq:I2wm311}
\end{eqnarray}
 We have checked that this series of reductions all converge for any
function \emph{f }that can be expanded in a powers series through
$\{n=-9,m=1,\nu=1\}$.

\section{Integral reductions for inverse powers in the exponentials }

\noindent Instead of the positive powers in the exponential of the
last section, let us take the reverse and set $p=q=h=j=0$,

\begin{eqnarray}
R_{2} & \left(\right. & n\left.,\,m,\,\nu,\,a,\,b,\,c,\,0,\,0,\,0,\,0\right)=\int_{0}^{\infty}\int_{0}^{\infty}\frac{1}{x^{n/2}y^{m/2}\left(x+y\right)^{\nu/2}}f\left(\frac{xy}{x+y}\right)\nonumber \\
 & \times & e^{-\frac{a}{x}-\frac{b}{y}-c\,xy/(x+y)}dx\,dy\quad.\label{eq:R2negpow}\\
 & = & \int_{0}^{\infty}\int_{0}^{\infty}X^{n/2}X^{m/2}\left(\frac{XY}{X+Y}\right)^{\nu/2}f\left(\frac{1}{X+Y}\right)\nonumber \\
 & \times & e^{-aX-bY-c/(X+Y)}dX\,dY\nonumber 
\end{eqnarray}
 where we have split off $e^{-c\,xy/(x+y)}$ from $f\left(\frac{xy}{x+y}\right)$
to explicate how this is different from simply making the replacement
$\left\{ x\to\frac{1}{X},\:y\to\frac{1}{Y},\:a\to p,\:b\to q\right\} $
in Eq. (\ref{eq:I2gen}). In Eq. (\ref{eq:I2stgen}) we apply a partial
fraction decomposition in the exponential of

\begin{equation}
-\frac{bs^{2}}{t\left(s^{2}+t\right)}=\frac{b}{s^{2}+t}-\frac{b}{t}\label{eq:splitb}
\end{equation}
 so that

\begin{eqnarray}
R_{2} & \left(\right. & n\left.,\,m,\,\nu,\,a,\,b,\,c,\,0,\,0,\,0,\,0\right)\nonumber \\
 & = & 2\int_{0}^{\infty}f(t)\,t^{-m/2}e^{-ct-qt-b/t}\,dt\int_{0}^{\infty}s^{m+\nu-3}\left(s^{2}+t\right)^{-\frac{m}{2}-\frac{n}{2}-\nu+2}\nonumber \\
 & \times & \exp\left(-\frac{a-b}{s^{2}+t}\right)ds\label{eq:R2negpows}
\end{eqnarray}

\noindent Changing variables to

\begin{equation}
w=\frac{a-b}{s^{2}+t}\label{eq:ssqtow}
\end{equation}
 transforms the \emph{s}-integral to

\begin{equation}
\frac{1}{a-b}\,\int_{0}^{\frac{a-b}{t}}e^{-w}\left(\frac{a-b}{w}\right)^{-\frac{m}{2}-\frac{n}{2}-\nu+4}\left(\frac{\sqrt{a-b-tw}}{\sqrt{w}}\right)^{m+\nu-4}dw\quad.\label{eq:J2negpow}
\end{equation}

The last term in the integrand is turned into a square if $m=\nu=3$
and the middle term is flipped upside-down if $n\geq1$ so letting
$\left\{ n=1,\:m=\nu=3\right\} $ gives \cite{GR5 p. 111 No. 2.311,GR5 p. 112 No. 2.322.1}


\begin{eqnarray}
R_{2} & \left(\right. & 1\left.,\,3,\,3,\,a,\,b,\,c,\,0,\,0,\,0,\,0\right)\nonumber \\
 & = & \int_{0}^{\infty}\int_{0}^{\infty}\frac{1}{x^{1/2}y^{3/2}\left(x+y\right)^{3/2}}f\left(\frac{xy}{x+y}\right)e^{-\frac{a}{x}-\frac{b}{y}-c\,xy/(x+y)}dx\,dy\label{eq:R2negpows133}\\
 & = & \left(a-b\right)^{-2}\int_{0}^{\infty}f(t)\,t^{-3/2}e^{-ct}\left(te^{-a/t}-e^{-b/t}(-a+b+t)\right)\,dt\nonumber 
\end{eqnarray}

\noindent and letting $\left\{ n=m=\nu=3\right\} $ gives \cite{GR5 p. 112 No. 2.322.1,GR5 p. 112 No. 2.322.2}

\begin{eqnarray}
R_{2} & \left(\right. & 3\left.3,\,3,\,3,\,a,\,b,\,c,\,0,\,0,\,0,\,0\right)\nonumber \\
 & = & \int_{0}^{\infty}\int_{0}^{\infty}\frac{1}{x^{3/2}y^{3/2}\left(x+y\right)^{3/2}}f\left(\frac{xy}{x+y}\right)e^{-\frac{a}{x}-\frac{b}{y}-c\,xy/(x+y)}dx\,dy\label{eq:R2negpows333}\\
 & = & \left(a-b\right)^{-3}\int_{0}^{\infty}f(t)\,t^{-3/2}e^{-ct}\left(e^{-b/t}(a-b-2t)+e^{-a/t}(a-b+2t)\right)\,dt\nonumber 
\end{eqnarray}
 and so on. The limit of what one may integrate is $e^{-w}\sqrt{\frac{1}{w}}(a-b-tw)$
so that $\left\{ n=0,\:m=\nu=3\right\} $ gives \cite{GR5 p. 364 No. 3.381.1}

\begin{eqnarray}
R_{2} & \left(\right. & 0\left.,\,3,\,3,\,a,\,b,\,c,\,0,\,0,\,0,\,0\right)\nonumber \\
 & = & \int_{0}^{\infty}\int_{0}^{\infty}\frac{1}{y^{3/2}\left(x+y\right)^{3/2}}f\left(\frac{xy}{x+y}\right)e^{-\frac{a}{x}-\frac{b}{y}-c\,xy/(x+y)}dx\,dy\label{eq:R2negpows033}\\
 & = & \left(a-b\right)^{-3/2}\int_{0}^{\infty}f(t)\,t^{-3/2}e^{-ct}\nonumber \\
 & \times & \left(\frac{1}{2}\sqrt{\pi}(2a-2b-t)e^{-b/t}\text{erf}\left(\sqrt{\frac{a-b}{t}}\right)+\sqrt{t}\sqrt{a-b}e^{-a/t}\right)\,dt\nonumber 
\end{eqnarray}

If instead we set the last term to unity with $m=\nu=2$ one gets
for $\left\{ n=1,\:m=\nu=2\right\} $ \cite{GR5 p. 362 No .3.361.1}

\begin{eqnarray}
R_{2} & \left(\right. & 1\left.,\,2,\,2,\,a,\,b,\,c,\,0,\,0,\,0,\,0\right)\nonumber \\
 & = & \int_{0}^{\infty}\int_{0}^{\infty}\frac{1}{x^{1/2}y\left(x+y\right)}f\left(\frac{xy}{x+y}\right)e^{-\frac{a}{x}-\frac{b}{y}-c\,xy/(x+y)}dx\,dy\nonumber \\
 & = & \left(a-b\right)^{-1/2}\int_{0}^{\infty}f(t)\,t^{-1}e^{-ct}e^{-b/t}\left(\sqrt{\pi}\text{erf}\left(\sqrt{\frac{a-b}{t}}\right)\right)\,dt\label{eq:R2negpows122}
\end{eqnarray}
 and for $\left\{ n=m=\nu=2\right\} $ \cite{GR5 p. 111 No. 2.311}

\vspace{0.1cm}

\begin{eqnarray}
R_{2} & \left(\right. & 1\left.,\,2,\,2,\,a,\,b,\,c,\,0,\,0,\,0,\,0\right)\nonumber \\
 & = & \int_{0}^{\infty}\int_{0}^{\infty}\frac{1}{xy\left(x+y\right)}f\left(\frac{xy}{x+y}\right)e^{-\frac{a}{x}-\frac{b}{y}-c\,xy/(x+y)}dx\,dy\nonumber \\
 & = & \left(a-b\right)^{-1}\int_{0}^{\infty}f(t)\,t^{-3/2}e^{-ct}\left(e^{-b/t}-e^{-a/t}\right)\,dt\quad.\label{eq:R2negpows222}
\end{eqnarray}
 The extension to larger values of $n$ is straightforward.

Finally, we note that none of the above are well defined for $a=b$.
Since the complicating exponential in the second line of (\ref{eq:R2negpows})
goes to unity in this case, we may write the general form

\begin{eqnarray}
R_{2} & \left(\right. & n\left.,\,m,\,\nu,\,b,\,b,\,c,\,0,\,0,\,0,\,0\right)\nonumber \\
 & = & \int_{0}^{\infty}\int_{0}^{\infty}\frac{1}{x^{n/2}y^{m/2}\left(x+y\right)^{\nu/2}}f\left(\frac{xy}{x+y}\right)e^{-\frac{b}{x}-\frac{b}{y}-c\,xy/(x+y)}dx\,dy\label{eq:R2negpows222-1}\\
 & = & \frac{\Gamma\left(\frac{1}{2}(m+\nu-2)\right)\Gamma\left(\frac{1}{2}(n+\nu-2)\right)}{\Gamma\left(\frac{1}{2}(m+n+2\nu-4)\right)}\int_{0}^{\infty}f(t)\,e^{-ct}e^{-b/t}t^{\frac{1}{2}(-m-n-\nu+2)}\,dt\;.\nonumber \\
 &  & \quad\Re(m+\nu)>2\land\Re(n+\nu)>2\nonumber 
\end{eqnarray}

\noindent This is most easily proved by setting $a=b$ and $h=0$
in Eq. (\ref{eq:R2p0q0m+nu>2_m+nu>2}) of section 5.

\section{Set of integral reductions for more complicated exponentials }

\noindent A entirely different set of integral reductions may be crafted
that include both positive and negative powers in the exponentials.
Using the same partial fraction decomposition as in Eq. (\ref{eq:splitb})
gives

\begin{eqnarray}
R{}_{2} & \left(\right. & n\left.,\,m,\,\nu,\,a,\,b,\,c,\,0,\,0,\,p,\,q\right)\nonumber \\
 & = & \int_{0}^{\infty}\int_{0}^{\infty}\frac{1}{y^{m/2}x^{n/2}\left(x+y\right)^{\nu/2}}f\left(\frac{xy}{x+y}\right)e^{-\frac{a}{x}-\frac{b}{y}-c\,xy/(x+y)-px-qy}dx\,dy\nonumber \\
 & = & 2\int_{0}^{\infty}f(t)\,t^{-m/2}e^{-ct-qt-b/t}\,dt\int_{0}^{\infty}s^{m+\nu-3}\left(s^{2}+t\right)^{-\frac{m}{2}-\frac{n}{2}-\nu+2}\nonumber \\
 & \times & \exp\left(-\frac{a-b}{s^{2}+t}-p\left(s^{2}+t\right)-\frac{qt^{2}}{s^{2}}\right)ds\label{eq:I2stgen-1}
\end{eqnarray}
 and we may well have convergence problems for $b>a$ unless \emph{c}
is large enough. Note that the exponential $e^{-c\,xy/(x+y)}$ multiplying
$f\left(\frac{xy}{x+y}\right)$ is transformed as $e^{-ct}$ multiplying
$f\left(t\right)$ as one would expect of a factor that can be folded
into the definition of $f$.

If we then specialize this integral to the case where $q=0$, one
may complete the square in the latter exponential and change variables

\vspace{0.04cm}

\begin{eqnarray}
\exp\left(-\frac{a-b}{s^{2}+t}-p\left(s^{2}+t\right)\right) & \hspace{-0.2cm}= & \hspace{-0.18cm}\exp\left(\hspace{-0.15cm}-\left(\sqrt{p}\sqrt{s^{2}+t}-\frac{\sqrt{a-b}}{\sqrt{s^{2}+t}}\right)^{2}\hspace{-0.2cm}-2\sqrt{p}\sqrt{a-b}\right)\nonumber \\
 & \equiv & \exp\left(-v^{2}-2\sqrt{p}\sqrt{a-b}\right)\quad.\label{eq:vsq}
\end{eqnarray}
 Of the four possible solutions, we choose the one with 
\[
s=-\sqrt{\frac{\sqrt{4\sqrt{p}v^{2}\sqrt{a-b}+4p(a-b)-4ap+4bp+v^{4}}}{2p}+\frac{\sqrt{a-b}}{\sqrt{p}}+\frac{v^{2}}{2p}-t}\qquad.
\]
 Then with

\begin{eqnarray}
ds & = & \frac{\left(s^{2}+t\right)^{3/2}}{s\left(\sqrt{a-b}+\sqrt{p}\left(s^{2}+t\right)\right)}dv\label{eq:ds2dv}
\end{eqnarray}
 we have \cite{Mathematica7}

\begin{eqnarray}
R_{2} & \left(\right. & n\left.,\,m,\,\nu,\,a,\,b,\,c,\,0,\,0,\,p,\,0\right)\nonumber \\
 & = & 2\int_{0}^{\infty}f(t)\,t^{-m/2}e^{-ct-b/t-2\sqrt{p}\sqrt{a-b}}\,dt\int_{\sqrt{p}\sqrt{t}-\sqrt{a-b}/\sqrt{t}}^{\infty}\:dv\,e^{-v^{2}}\nonumber \\
 & \times & \frac{\sqrt{p}2^{\frac{1}{2}(n+\nu-1)}\left(\frac{\sqrt{4\sqrt{p}v^{2}\sqrt{a-b}+v^{4}}+2\sqrt{p}\sqrt{a-b}+v^{2}}{p}\right)^{\frac{1}{2}(-m-n-2\nu+7)}}{\left(\sqrt{4\sqrt{p}v^{2}\sqrt{a-b}+v^{4}}+4\sqrt{p}\sqrt{a-b}+v^{2}\right)}\nonumber \\
 & \times & \frac{1}{\left(-\sqrt{\frac{\sqrt{4\sqrt{p}v^{2}\sqrt{a-b}+v^{4}}+2\sqrt{p}\sqrt{a-b}-2pt+v^{2}}{p}}\right)^{4-m-\nu}}\quad.\label{eq:I2q0}
\end{eqnarray}

If we set $\{m=4-\nu\}$, the last factor of (\ref{eq:I2q0}) goes
to unity, giving us some hope of doing the integral. But even the
simplest of those, with $\{n=3-\nu\}$ \cite{Mathematica7}

\vspace{4.5cm}

\begin{eqnarray}
R_{2} & \left(\right. & 3-\nu\left.,\,4-\nu,\,\nu,\,a,\,b,\,c,\,0,\,0,\,p,\,0\right)\nonumber \\
 & = & 2\int_{0}^{\infty}f(t)\,t^{-\left(4-\nu\right)/2}e^{-ct-b/t-2\sqrt{p}\sqrt{a-b}}\,dt\int_{\sqrt{p}\sqrt{t}-\sqrt{a-b}/\sqrt{t}}^{\infty}\:dv\,e^{-v^{2}}\nonumber \\
 & \times & \frac{\sqrt{p}2}{\left(\sqrt{4\sqrt{p}v^{2}\sqrt{a-b}+v^{4}}+4\sqrt{p}\sqrt{a-b}+v^{2}\right)}\label{eq:R2q0m4-nu_n3-nu}\\
 & = & 2\int_{0}^{\infty}f(t)\,t^{-m/2}e^{-ct-b/t-2\sqrt{p}\sqrt{a-b}}\,dt\nonumber \\
 & \times & \left(\frac{t\sqrt{a-b}e^{4\sqrt{p}\sqrt{a-b}}\sqrt{\frac{\left(-a+b+pt^{2}\right)^{2}}{t^{2}}}\text{erf}\left(\frac{\sqrt{2\sqrt{p}t\sqrt{a-b}+a-b+pt^{2}}}{\sqrt{t}}\right)}{\sqrt{2\sqrt{p}t\sqrt{a-b}+a-b+pt^{2}}}\right.\nonumber \\
 & + & \left.\left(\sqrt{p}t\sqrt{a-b}-a+b\right)\text{erf}\left(\frac{\sqrt{a-b}-\sqrt{p}t}{\sqrt{t}}\right)\right.\nonumber \\
 & + & \left.\left(e^{4\sqrt{p}\sqrt{a-b}}-1\right)\left(-\sqrt{p}t\sqrt{a-b}+a-b\right)\right)\nonumber \\
 & \times & \frac{\sqrt{\pi}}{8\sqrt{p}\sqrt{a-b}\left(\sqrt{p}t\sqrt{a-b}-a+b\right)}\nonumber \\
 &  & {\scriptstyle {\left[\left(\frac{\sqrt{a-b}-\sqrt{p}t}{\sqrt{t}}\notin\mathbb{R}\lor\Re\left(\frac{\sqrt{a-b}-\sqrt{p}t}{\sqrt{t}}\right)\leq0\right),\right.}}\nonumber \\
 &  & {\scriptstyle {\left(\frac{-2\sqrt{\sqrt{p}t^{2}\left(-\sqrt{a-b}\right)}+\sqrt{t}\sqrt{a-b}-\sqrt{p}t^{3/2}}{t}\notin\mathbb{R}\lor\Re\left(\frac{-2\sqrt{\sqrt{p}t^{2}\left(-\sqrt{a-b}\right)}+\sqrt{t}\sqrt{a-b}-\sqrt{p}t^{3/2}}{t}\right)\leq0\right)}}\nonumber \\
 &  & {\scriptstyle {\left.\left(\frac{2\sqrt{\sqrt{p}t^{2}\left(-\sqrt{a-b}\right)}+\sqrt{t}\sqrt{a-b}-\sqrt{p}t^{3/2}}{t}\notin\mathbb{R}\lor\Re\left(\frac{2\sqrt{\sqrt{p}t^{2}\left(-\sqrt{a-b}\right)}+\sqrt{t}\sqrt{a-b}-\sqrt{p}t^{3/2}}{t}\right)\leq0\right)\right]}}\nonumber 
\end{eqnarray}
 has conditions that cannot be met as $t$ becomes very small unless$\sqrt{p}t=\sqrt{a-b}$
or one of these coefficients in the original exponents is an imaginary
number.

So let us return to the integral in terms of \emph{x} and $y$ to
form a self-consistent version of the former condition by setting
$p=(a-b)(x+y)^{2}/\left(x^{2}y^{2}\right)$. The utility of this value
is not apparent in terms of the original variables (where $-\frac{a}{x}-\frac{b}{y}-px\to\frac{x(b-a)}{y^{2}}+\frac{b-2a}{x}+\frac{b-2a}{y}$)
but creates a useful exponential when expressed in terms of \emph{s}
and \emph{t} on the last line of the modified integral

\vspace{2.4cm}

\begin{eqnarray}
\tilde{R}{}_{2} & \left(\right. & n\left.,\,m,\,\nu,\,a,\,b,\,c,\,q\right)=\hspace{-0cm}\hspace{-0cm}\int_{0}^{\infty}\hspace{-0.2cm}\int_{0}^{\infty}\hspace{-0cm}\frac{1}{y^{m/2}x^{n/2}\left(x+y\right)^{\nu/2}}f\left(\frac{xy}{x+y}\right)\nonumber \\
 & \times & e^{-a/x-b/y-cxy/(x+y)-(a-b)(x+y)^{2}/(xy^{2})-qy}dx\,dy\nonumber \\
 & = & 2\int_{0}^{\infty}f(t)\,t^{-m/2}e^{-ct-qt-b/t}\,dt\int_{0}^{\infty}s^{m+\nu-3}\left(s^{2}+t\right)^{-\frac{m}{2}-\frac{n}{2}-\nu+2}\nonumber \\
 & \times & \exp\left(-\frac{a-b}{s^{2}+t}-\frac{(a-b)\left(s^{2}+t\right)}{t^{2}}-\frac{qt^{2}}{s^{2}}\right)ds\label{eq:I2stgen-1-1}
\end{eqnarray}

We again specialize this integral to the case where $q=0$ and complete
the square in the latter exponential and change variables

\begin{eqnarray}
{\textstyle \exp\left(-\frac{a-b}{s^{2}+t}-\frac{(a-b)\left(s^{2}+t\right)}{t^{2}}\right)} & {\textstyle =} & {\textstyle \exp\left(-\left(\frac{\sqrt{a-b}\sqrt{s^{2}+t}}{t}-\frac{\sqrt{a-b}}{\sqrt{s^{2}+t}}\right)^{2}-\frac{2(a-b)}{t}\right)}\nonumber \\
 & {\textstyle \equiv} & {\textstyle \exp\left(-v^{2}-\frac{2(a-b)}{t}\right)}\quad.\label{eq:vsq-1}
\end{eqnarray}

\noindent Of the four possible solutions, we choose the one with 
\begin{equation}
s=-\sqrt{\frac{t^{3/2}v\sqrt{4a-4b+tv^{2}}}{2(a-b)}+\frac{t^{2}v^{2}}{2(a-b)}}\qquad.\label{eq:stov2}
\end{equation}

\noindent Then with

\begin{eqnarray}
\text{ds} & = & \frac{s\sqrt{a-b}\left(s^{2}+2t\right)}{t\left(s^{2}+t\right)^{3/2}}\text{dv}\label{eq:ds2dv-1}
\end{eqnarray}
 we have \cite{Mathematica7}

\begin{eqnarray}
\tilde{R}_{2} & \left(\right. & n\left.,\,m,\,\nu,\,a,\,b,\,c,\,p,\,0\right)\label{eq:Ralt2q0}\\
 & \hspace{-0.8cm}= & \hspace{-0.4cm}2\int_{0}^{\infty}f(t)\,t^{-m/2}e^{-ct-b/t-2(a-b)/t}\,dt\int_{0}^{\infty}\:dv\,e^{-v^{2}}\left(\frac{1}{a-b}\right)^{\frac{1}{2}(-n-\nu+2)}\nonumber \\
 & \hspace{-0.8cm}\times & \hspace{-0.4cm}\frac{(-1)^{m+\nu-2}2^{\frac{1}{2}(n+\nu-1)}\left(t\left(\sqrt{t}v\sqrt{4a-4b+tv^{2}}+2a-2b+tv^{2}\right)\right)^{\frac{1}{2}(-m-n-2\nu+7)}}{\left(\sqrt{t}v\sqrt{4a-4b+tv^{2}}+4a-4b+tv^{2}\right)\left(t^{3/2}v\sqrt{4a-4b+tv^{2}}+t^{2}v^{2}\right)^{\frac{1}{2}(4-m-\nu)}}\nonumber 
\end{eqnarray}
 If we set $\{m=4-\nu\}$ the second factor in the denominator of
(\ref{eq:Ralt2q0}) goes to unity, as does the last factor in the
numerator if we then set $\{n=3-\nu\}$, giving \cite{Mathematica7}

\vspace{1.3cm}

\begin{eqnarray}
\tilde{R}{}_{2} & \left(\right. & 3-\nu\left.,\,4-\nu,\,\nu,\,a,\,b,\,c,\,q\right)=\int_{0}^{\infty}\int_{0}^{\infty}\frac{1}{x^{\left(3-\nu\right)/2}y^{\left(4-\nu\right)/2}\left(x+y\right)^{\nu/2}}f\left(\frac{xy}{x+y}\right)\nonumber \\
 & \hspace{-0.8cm}\times & \hspace{-0.4cm}e^{-a/x-b/y-cxy/(x+y)-(a-b)(x+y)^{2}/(xy^{2})-qy}dx\,dy\label{eq:Raltn3-nu_m4-nu}\\
 & \hspace{-0cm}= & \hspace{-0.1cm}2\int_{0}^{\infty}f(t)\,t^{-m/2}e^{-ct-b/t-2(a-b)/t}\frac{\sqrt{\pi}}{4\sqrt{a-b}}\left(1-e^{4(a-b)/t}\text{erfc}\left(\frac{2\sqrt{a-b}}{\sqrt{t}}\right)\right)\,dt\quad.\nonumber \\
 &  & {\textstyle {\quad\left[\sqrt{\frac{\Im(b)-\Im(a)}{\Im(t)}}\le0\lor\Re(a)\Im(t)+\Im(b)\Re(t)\geq\Im(a)\Re(t)+\Re(b)\Im(t),\right.}}\nonumber \\
 &  & {\textstyle {\quad\frac{\sqrt{\Im(b)-\Im(a)}}{\sqrt{\ Im(t)}}\geq0\lor\Re(a)+\frac{\Im(b)\Re(t)}{\Im(t)}\geq\frac{\Im(a)\Re(t)}{\Im(t)}+\Re(b),}}\nonumber \\
 &  & {\textstyle {\quad\left.\frac{\sqrt{b-a}}{\ sqrt{t}}\notin\mathbb{R}\lor\left(\frac{\sqrt{b-a}}{\sqrt{t}}\neq0\land\Re\left(\frac{\sqrt{b-a}}{\sqrt{t}}\right)=0\right)\right]}}\nonumber 
\end{eqnarray}
 The powers of $\nu$ can be folded into a redefinition of \emph{f}.

If we instead retain the last factor in the numerator with unit power
by setting $\{n=1-\nu\}$, after reducing we obtain \cite{Mathematica7}

\begin{eqnarray}
\tilde{R}{}_{2} & \left(\right. & 1-\nu\left.,\,4-\nu,\,\nu,\,a,\,b,\,c,\,q\right)\nonumber \\
 & = & \int_{0}^{\infty}\int_{0}^{\infty}\frac{1}{x^{\left(1-\nu\right)/2}y^{\left(4-\nu\right)/2}\left(x+y\right)^{\nu/2}}f\left(\frac{xy}{x+y}\right)\nonumber \\
 & \times & e^{-a/x-b/y-cxy/(x+y)-(a-b)(x+y)^{2}/(xy^{2})-qy}dx\,dy\nonumber \\
 & = & 2\int_{0}^{\infty}f(t)\,\frac{\sqrt{\pi}t^{\frac{\nu}{2}-1}}{2\sqrt{a-b}}e^{-2(a-b)/t-5b/t-ct}\nonumber \\
 & \times & \left(-e^{4a/t}\text{erf}\left(\frac{2\sqrt{a-b}}{\sqrt{t}}\right)+e^{4a/t}+e^{4b/t}\right)\,dt\quad.\label{eq:Raltn1-nu_m4-nu}\\
 &  & {\scriptstyle \quad\left[\Re(a)\geq\Re(b),\Re(a)+\frac{\Im(b)\Re(t)}{\Im(t)}\geq\frac{\Im(a)\Re(t)}{\Im(t)}+\Re(b)\right.}\nonumber \\
 &  & {\scriptstyle \left.\lor\right.\hspace{-0.1cm}{\scriptstyle {\left.\left(\left(\frac{\sqrt{\Im(b)-\Im(a)}}{\sqrt{\Im(t)}}\notin\mathbb{R}\lor\Re\left(\frac{\sqrt{\Im(b)-\Im(a)}}{\sqrt{\Im(t)}}\right)=0\right)\land\left(\sqrt{\frac{\Im(b)-\Im(a)}{\Im(t)}}\notin\mathbb{R}\lor\Re\left(\sqrt{\frac{\Im(b)-\Im(a)}{\Im(t)}}\right)\leq0\right)\right)\right]}}}\nonumber 
\end{eqnarray}

One can also move the second factor in the denominator of (\ref{eq:Ralt2q0})
into the numerator with unit power by setting $\{m=6-\nu\}$. We then
set $\{n=1-\nu\}$ to eliminate the other term in the numerator, and
reduce the resultant quotient to obtain \cite{Mathematica7}

\begin{eqnarray}
\tilde{R}{}_{2} & \left(\right. & 1-\nu\left.,\,6-\nu,\,\nu,\,a,\,b,\,c,\,q\right)\nonumber \\
 & = & \int_{0}^{\infty}\int_{0}^{\infty}\frac{1}{x^{\left(1-\nu\right)/2}y^{\left(6-\nu\right)/2}\left(x+y\right)^{\nu/2}}f\left(\frac{xy}{x+y}\right)\nonumber \\
 & \times & e^{-a/x-b/y-cxy/(x+y)-(a-b)(x+y)^{2}/(xy^{2})-qy}dx\,dy\nonumber \\
 & = & 2\int_{0}^{\infty}f(t)\,t^{-2+\nu/2}e^{-ct-qt-b/t-2(a-b)/t}\frac{\sqrt{\pi}}{2\sqrt{a-b}}\text{erfc}\left(\frac{2\sqrt{a-b}}{\sqrt{t}}\right)\,dt\quad.\label{eq:Raltn1-nu_m6-nu}\\
 &  & {\scriptstyle \quad\left[\Re(a-b)>0,\right.}\nonumber \\
 &  & {\scriptstyle \quad\Re(a)+\frac{\Im(b)\Re(t)}{\Im(t)}\geq\frac{\Im(a)\Re(t)}{\Im(t)}+\Re(b)\lor\left(\sqrt{\frac{\Im(b)-\Im(a)}{\Im(t)}}\leq0\land\frac{\sqrt{\Im(b)-\Im(a)}}{\sqrt{\Im(t)}}\geq0\right),}\nonumber \\
 &  & {\scriptstyle \quad\left.\frac{a-b}{t}\notin\mathbb{R}\lor\Re\left(\frac{a-b}{t}\right)\geq0\right]}\nonumber 
\end{eqnarray}
 Again, the powers of $\nu$ can be folded into a redefinition of
\emph{f}.

One may also move the second factor in the denominator of (\ref{eq:Ralt2q0})
into the numerator with a power of two by setting $\{m=8-\nu\}$.
Again eliminate the other term in the numerator with $\{n=1-\nu\}$,
and reduce the resultant quotient somewhat. Mathematica 7 was unable
to do this integral but Mathematica 9 could, with

\vspace{-0.6cm}

\begin{eqnarray}
\tilde{R}{}_{2} & \left(\right. & 1-\nu\left.,\,8-\nu,\,\nu,\,a,\,b,\,c,\,q\right)\nonumber \\
 & = & \int_{0}^{\infty}\int_{0}^{\infty}\frac{1}{x^{\left(1-\nu\right)/2}y^{\left(8-\nu\right)/2}\left(x+y\right)^{\nu/2}}f\left(\frac{xy}{x+y}\right)\nonumber \\
 & \times & e^{-a/x-b/y-cxy/(x+y)-(a-b)(x+y)^{2}/(xy^{2})-qy}dx\,dy\nonumber \\
 & = & 2\int_{0}^{\infty}f(t)\,t^{-2+\nu/2}e^{-ct-qt-b/t-2(a-b)/t}\,dt\frac{t^{\frac{\nu}{2}-3}}{8(a-b)^{2}}\nonumber \\
 & \times & \left(\sqrt{\pi}\sqrt{a-b}e^{-\frac{4b}{t}}\left(-8be^{\frac{4a}{t}}+t\left(e^{\frac{4b}{t}}-e^{\frac{4a}{t}}\right)+8ae^{\frac{4a}{t}}\right)\right.\quad.\label{eq:Raltn1-nu_m8-nu}\\
 & - & \left.e^{-\frac{4b}{t}}\left(\sqrt{\pi}\sqrt{a-b}e^{\frac{4a}{t}}(8a-8b-t)\text{erf}\left(\frac{2\sqrt{a-b}}{\sqrt{t}}\right)+4\sqrt{t}(a-b)e^{\frac{4b}{t}}\right)\right)\nonumber \\
 &  & {\scriptstyle \quad\left[\Re(a)\geq\Re(b),\Re(a)+\frac{\Im(b)\Re(t)}{\Im(t)}\geq\frac{\Im(a)\Re(t)}{\Im(t)}+\Re(b)\right.}\nonumber \\
 &  & {\scriptstyle \left.\lor\right.}\hspace{-0.1cm}{\scriptstyle {\left.\left(\left(\frac{\sqrt{\Im(b)-\Im(a)}}{\sqrt{\Im(t)}}\notin\mathbb{R}\lor\Re\left(\frac{\sqrt{\Im(b)-\Im(a)}}{\sqrt{\Im(t)}}\right)=0\right)\land\left(\sqrt{\frac{\Im(b)-\Im(a)}{\Im(t)}}\notin\mathbb{R}\lor\Re\left(\sqrt{\frac{\Im(b)-\Im(a)}{\Im(t)}}\right)\leq0\right)\right)\right]}}\nonumber 
\end{eqnarray}

Finally, one can also move the second factor in the denominator of
(\ref{eq:Ralt2q0}) into the numerator with unit power by setting
$\{m=6-\nu\}$ while allowing the other term in the numerator to appear
with unit power by setting $\{n=-1-\nu\}$, and reduce the resultant
quotient, one obtains

\vspace{-0.6cm}

\begin{eqnarray}
\tilde{R}{}_{2}\hspace{-0cm} & \left(\right. & -1-\nu\left.,\,6-\nu,\,\nu,\,a,\,b,\,c,\,q\right)\nonumber \\
 & = & \int_{0}^{\infty}\int_{0}^{\infty}\frac{1}{x^{\left(-1-\nu\right)/2}y^{\left(6-\nu\right)/2}\left(x+y\right)^{\nu/2}}f\left(\frac{xy}{x+y}\right)\nonumber \\
 & \times & e^{-a/x-b/y-cxy/(x+y)-(a-b)(x+y)^{2}/(xy^{2})-qy}dx\,dy\nonumber \\
 & = & \hspace{-0cm}2\int_{0}^{\infty}f(t)\,t^{-1+\nu/2}e^{-ct-qt-b/t-2(a-b)/t}\,dt\nonumber \\
 & \times & \left(e^{-\frac{4b}{t}}\left(\sqrt{\pi}e^{\frac{4a}{t}}(4a-4b-t)\text{erf}\left(\frac{2\sqrt{a-b}}{\sqrt{t}}\right)+4\sqrt{t}\sqrt{a-b}e^{\frac{4b}{t}}\right)\right.\nonumber \\
 & + & \left.\left(\sqrt{\pi}e^{-\frac{4b}{t}}\left(4be^{\frac{4a}{t}}+t\left(e^{\frac{4a}{t}}+e^{\frac{4b}{t}}\right)-4ae^{\frac{4a}{t}}\right)\right)\frac{1}{4(a-b)^{3/2}}\right)\quad.\label{eq:Raltn-1-nu_m6-nu-2}\\
 &  & {\scriptstyle \quad\left[\Re(a)\geq\Re(b),\Re(a)+\frac{\Im(b)\Re(t)}{\Im(t)}\geq\frac{\Im(a)\Re(t)}{\Im(t)}+\Re(b)\lor\right.}\nonumber \\
 &  & \hspace{-0.1cm}{\scriptstyle {\left.\left(\left(\frac{\sqrt{\Im(b)-\Im(a)}}{\sqrt{\Im(t)}}\notin\mathbb{R}\lor\Re\left(\frac{\sqrt{\Im(b)-\Im(a)}}{\sqrt{\Im(t)}}\right)=0\right)\land\left(\sqrt{\frac{\Im(b)-\Im(a)}{\Im(t)}}\ notin\mathbb{R}\lor\Re\left(\sqrt{\frac{\Im(b)-\Im(a)}{\Im(t)}}\right)\leq0\right)\right)\right]}}\nonumber 
\end{eqnarray}

\section{Set of integral reductions for inverse powers and powers times inverse
binomials in the exponentials }

\subsection{Integrals with $j=p=q=0$ }

\noindent An entirely different set of integral reductions may be
crafted for inverse powers in the exponentials. If we set $j=p=q=0$
in (\ref{eq:I2stgen}) we obtain

\begin{eqnarray}
R{}_{2} & \left(\right. & n\left.,\,m,\,\nu,\,a,\,b,\,c,\,h,\,0,\,0,\,0\right)\hspace{-0.3cm}\nonumber \\
 & = & \hspace{-0.1cm}\int_{0}^{\infty}\int_{0}^{\infty}\frac{1}{y^{m/2}x^{n/2}\left(x+y\right)^{\nu/2}}f\left(\frac{xy}{x+y}\right)e^{-\frac{a}{x}-\frac{b}{y}-c\,xy/(x+y)-h\,y/(x+y)}dx\,dy\nonumber \\
 & = & 2\int_{0}^{\infty}f(t)\,t^{-m/2}e^{-ct-b/t}\,dt\int_{0}^{\infty}s^{m+\nu-3}\left(s^{2}+t\right)^{-\frac{m}{2}-\frac{n}{2}-\nu+2}\nonumber \\
 & \times & \exp\left(-\frac{a-b+ht}{s^{2}+t}\right)ds\quad.\label{eq:R2stp0q0}
\end{eqnarray}
 If we change variables to

\begin{equation}
w=\frac{a-b+ht}{s^{2}+t},\:s=\frac{\sqrt{a-b+ht-tw}}{\sqrt{w}}\qquad,\label{eq:stow}
\end{equation}
 with

\begin{eqnarray}
ds & = & -\frac{\left(s^{2}+t\right)^{2}}{2s(a-b+ht)}dw\label{eq:ds2dv-1-1}
\end{eqnarray}
 we have

\begin{eqnarray}
R{}_{2} & \left(\right. & n\left.,\,m,\,\nu,\,a,\,b,\,c,\,h,\,0,\,0,\,0\right)\hspace{-0.3cm}\nonumber \\
 & = & \hspace{-0.1cm}\int_{0}^{\infty}\int_{0}^{\infty}\frac{1}{x^{n/2}y^{m/2}\left(x+y\right)^{\nu/2}}f\left(\frac{xy}{x+y}\right)e^{-\frac{a}{x}-\frac{b}{y}-c\,xy/(x+y)-h\,y/(x+y)}dx\,dy\nonumber \\
 & = & 2\int_{0}^{\infty}f(t)\,t^{-m/2}e^{-ct-b/t}(a-b+ht)^{-m/2-n/2-\nu+3}\,dt\nonumber \\
 & \times & \int_{0}^{\left(a-b+ht\right)/2}w^{n/2+\nu/2-2}e^{-w}(a-b+ht-tw)^{m/2+\nu/2-2}dw\label{eq:R2p0q0w}
\end{eqnarray}

The general result, good for $(m+\nu)>2,\:(n+\nu)>2$ is

\vspace{2.7cm}

\begin{eqnarray}
R{}_{2} & \left(\right. & n\left.,\,m,\,\nu,\,a,\,b,\,c,\,h,\,0,\,0,\,0\right)\hspace{-0.3cm}\nonumber \\
 & = & \hspace{-0.1cm}\int_{0}^{\infty}\int_{0}^{\infty}\frac{1}{x^{n/2}y^{m/2}\left(x+y\right)^{\nu/2}}f\left(\frac{xy}{x+y}\right)e^{-\frac{a}{x}-\frac{b}{y}-c\,xy/(x+y)-h\,y/(x+y)}dx\,dy\nonumber \\
 & = & \int_{0}^{\infty}f(t)\,e^{-\frac{b}{t}-ct}t^{\frac{1}{2}(-m-n-\nu+2)}\frac{\Gamma\left(\frac{1}{2}(m+\nu-2)\right)\Gamma\left(\frac{1}{2}(n+\nu-2)\right)}{\Gamma\left(\frac{1}{2}(m+n+2\nu-4)\right)}\label{eq:R2p0q0m+nu>2_m+nu>2}\\
 &  & \,_{1}F_{1}\left(\frac{1}{2}(n+\nu-2);\frac{1}{2}(m+n+2\nu-4);-\frac{a-b+ht}{t}\right)\,dt\nonumber \\
 &  & \quad\left[\Re(m+\nu)>2,\,\Re(n+\nu)>2,\{a,\,b,\,t,\,h\}\in\mathbb{R},\,t>0,\right.\nonumber \\
 &  & \quad\left.((a=b\land h>0)\lor((a-b+th)>0\land b\neq a)\right]\nonumber 
\end{eqnarray}
 Whenever $m=n\pm2M$, where M is a non-negative integer, we may simplify
this with \cite{PBM3 p. 579 No. 7.11.5-7}

\begin{eqnarray}
\,_{1}F_{1} & ( & A;2A-M;z)=\Gamma\left(A-M-\frac{1}{2}\right)\left(\frac{z}{4}\right)^{M-A+\frac{1}{2}}e^{z/2}\label{eq:1F1toI-M}\\
 & \hspace{-0.3cm}\times & \sum_{k=0}^{M}\frac{\left((-1)^{k}(-M)_{k}(2A-2M-1)_{k}\right)\left(A+k-M-\frac{1}{2}\right)I_{A+k-M-\frac{1}{2}}\left(\frac{z}{2}\right)}{(2A-M)_{k}k!}\quad,\nonumber 
\end{eqnarray}

\begin{equation}
\,_{1}F_{1}(A;2A;z)=2^{2A-1}e^{z/2}(-z)^{\frac{1}{2}-A}\Gamma\left(A+\frac{1}{2}\right)I_{A-\frac{1}{2}}\left(-\frac{z}{2}\right)\quad,\label{eq:1F1toI}
\end{equation}
 and

\begin{eqnarray}
\,_{1}F_{1}(A;2A+M;z) & = & \Gamma\left(A-\frac{1}{2}\right)\left(\frac{z}{4}\right)^{-A+\frac{1}{2}}e^{z/2}\nonumber \\
 & \times & \sum_{k=0}^{M}\frac{\left((-M)_{k}(2A-1)_{k}\right)\left(A+k-\frac{1}{2}\right)I_{A+k-\frac{1}{2}}\left(\frac{z}{2}\right)}{(2A+M)_{k}k!}\quad,\label{eq:1F1toI+M}
\end{eqnarray}
 that can be readily integrated whenever $n+\nu=2\left(A+1\right)$
is an even integer. Even when $\nu$ alone is an even integer, and
M is a non-negative integer, we may nevertheless find an integrable
form with

\begin{equation}
\,_{1}F_{1}(A;A-M;z)=\frac{(-1)^{M}e^{z}M!L_{M}^{A-M-1}(-z)}{(1-A)_{M}}\quad.\label{eq:1F1toL}
\end{equation}

\subsection{Application}

\noindent Probability amplitudes in atomic physics and variational
wave functions for many-electron atoms involve products of Yukawa
or (their special case) Coulomb potentials and hydrogenic orbitals
that may be derived from Yukawa potentials via derivatives. Consider,
then, the case of two Yukawa potentials centered on different positions
\begin{equation}
S_{1}^{\eta_{1}0\eta_{2}0}\left(0,;0,x_{2}\right)=\int d^{3}x_{1}\frac{e^{-\eta_{1}x_{1}}}{x_{1}}\frac{e^{-\eta_{2}x_{12}}}{x_{12}}\equiv\int d^{3}x_{1}\frac{e^{-\eta_{1}x_{1}}}{x_{1}}\frac{e^{-\eta_{2}\left|\mathbf{x}_{1}-\mathbf{x}_{2}\right|}}{\left|\mathbf{x}_{1}-\mathbf{x}_{2}\right|}\label{eq:SYYx2}
\end{equation}

The author has given an analytically reduced form for multidimensional
integrals over any number of such products in terms of Gaussian transforms,
generally with one remaining integral for each atomic center in the
original integrals. \cite{Stra89a} Such a reduction in the present
case from three to two integral dimensions is in no way dramatic,
but does provide a soluble example for the utility of the (further)
reduction formulae of this section. One may write down the final form
for the above integral using the notation formalism in that paper,
but the reduction may instead be easily found by completing the square
in the Gaussian transform

\begin{eqnarray}
S_{1}^{\eta_{1}0\eta_{2}0} & \left(\right. & 0\left.,;0,x_{2}\right)\nonumber \\
 & = & \int d^{3}x_{1}{\displaystyle \frac{1}{\sqrt{\pi}}}\int_{0}^{\infty}\,d\rho_{1}{\displaystyle \frac{e^{-x_{1}^{2}\rho_{1}}e^{-\eta_{1}^{2}/4/\rho_{1}}}{\rho_{1}^{\;1/2}}}{\displaystyle \frac{1}{\sqrt{\pi}}}\int_{0}^{\infty}\,d\rho_{2}{\displaystyle \frac{e^{-x_{12}^{2}\rho_{2}}e^{-\eta_{2}^{2}/4/\rho_{2}}}{\rho_{2}^{\;1/2}}}\nonumber \\
 & = & {\displaystyle \frac{1}{\pi}}\int d^{3}x'_{1}\int_{0}^{\infty}\,d\rho_{1}{\displaystyle \frac{e^{-\eta_{1}^{2}/4/\rho_{1}}}{\rho_{1}^{\;1/2}}}\int_{0}^{\infty}\,d\rho_{2}{\displaystyle \frac{e^{-\eta_{2}^{2}/4/\rho_{2}}}{\rho_{2}^{\;1/2}}}\nonumber \\
 & \times & exp\left(-\left(\rho_{1}+\rho_{2}\right)x'{}_{1}^{2}-\frac{x_{2}^{2}\rho_{1}\rho_{2}}{\rho_{1}+\rho_{2}}\right)\label{eq:SYYx2Gaussian}
\end{eqnarray}
 where we have changed variables from $\mathbf{x}_{1}$ to $\mathbf{x}'_{1}=\mathbf{x}_{1}-\frac{\rho_{2}}{\rho_{1}+\rho_{2}}\mathbf{x}_{2}$
with unit Jacobian. Then the spatial integral may be done \cite{GR5 p. 382 No. 3.461.2}

\begin{eqnarray}
\int e^{-\left(\rho_{1}+\rho_{2}\right)x  '^2 _1}d^{3} x '_{1}&=& 4\pi\int_{0}^{\infty}e^{-\left(\rho_{1}+\rho_{2}\right)x '^2 _1}\; x  '^2 _1 dx ' _1=\frac{4\pi^{1+1/2}}{2^{2}\left(\rho_{1}+\rho_{2}\right)^{3/2}}\quad. \\
 &  &\left[\rho_{1}+\rho_{2}>0\right]  \nonumber
 \end{eqnarray}

What is left is

\begin{eqnarray}
S_{1}^{\eta_{1}0\eta_{2}0} & \left(\right. & 0\left.,;0,x_{2}\right)=\pi^{1/2}\int_{0}^{\infty}\,d\rho_{1}{\displaystyle \frac{e^{-\eta_{1}^{2}/4/\rho_{1}}}{\rho_{1}^{\;1/2}}}\int_{0}^{\infty}\,d\rho_{2}{\displaystyle \frac{e^{-\eta_{2}^{2}/4/\rho_{2}}}{\rho_{2}^{\;1/2}}}\label{eq:threerphos-1-1}\\
 & \times & \frac{1}{\left(\rho_{1}+\rho_{2}\right)^{3/2}}\exp\left(-\frac{x_{2}^{2}\rho_{1}\rho_{2}}{\rho_{1}+\rho_{2}}\right)\nonumber \\
 & =\pi^{1/2} & R{}_{2}\left(4,\,4,\,0,\,\eta_{1}^{2}/4,\,\eta_{2}^{2}/4,\,x_{2}^{2},\,0,\,0,\,0,\,0\right)\nonumber 
\end{eqnarray}
 where in the last line we have taken our function to be 
\begin{equation}
f\left(t\right)=t^{\mu}\label{eq:fispowers}
\end{equation}
 and in the present case have $\mu=3/2$. We could also have set $\mu=0$
with $n=m=1$ and $\nu=3$. Then using (\ref{eq:1F1toI}) we have

\begin{eqnarray}
\pi^{1/2}R{}_{2} & \left(\right. & 4\left.,\,4,\,0,\,\eta_{1}^{2}/4,\,\eta_{2}^{2}/4,\,x_{2}^{2},\,0,\,0,\,0\right)\nonumber \\
 & = & \pi^{1/2}\int_{0}^{\infty}dt\left(\frac{4e^{-tx_{2}^{2}}}{\sqrt{t}\left(\eta_{2}^{2}-\eta_{1}^{2}\right)}\right)\left(e^{-\eta_{1}^{2}/4/t}-e^{-\eta_{2}^{2}/4/t}\right)\nonumber \\
 & = & \frac{4\pi\left(e^{-x_{2}\eta_{1}}-e^{-x_{2}\eta_{2}}\right)}{x_{2}\left(\eta_{2}^{2}-\eta_{1}^{2}\right)}\quad,\label{eq:R2p0q0m+nu>2_m+nu>2-1}\\
 &  & \quad\left[\Re\left(x_{2}^{2}\right)\geq0,\,\Re\left(\eta_{1}^{2}\right)\geq0,\,\Re\left(\eta_{2}^{2}\right)\geq0\right]\nonumber 
\end{eqnarray}
 which is indeed the correct result.

Suppose instead of a Yukawa potential in $\eta_{1}$ we have a hydrogenic
\emph{1s} wave function, that may be had by differentiation:

\begin{eqnarray}
S_{1s,1}^{\eta_{1}\eta_{2}0}\left(0,;0,x_{2}\right) & = & \int d^{3}x_{1}u_{1s}^{\eta_{1}}\left(x_{1}\right)\frac{e^{-\eta_{2}x_{12}}}{x_{12}}\nonumber \\
 & = & \int d^{3}x_{1}\frac{\eta_{1}^{3/2}e^{-\eta_{1}x_{1}}}{\sqrt{\pi}}\frac{e^{-\eta_{2}x_{12}}}{x_{12}}\nonumber \\
 & = & -\frac{\eta_{1}^{3/2}}{\sqrt{\pi}}\frac{\partial}{\partial\eta_{1}}\int d^{3}x_{1}\frac{e^{-\eta_{1}x_{1}}}{x_{1}}\frac{e^{-\eta_{2}x_{12}}}{x_{12}}\nonumber \\
 & = & -\frac{\eta_{1}^{3/2}}{\sqrt{\pi}}\frac{\partial}{\partial\eta_{1}}S_{1}^{\eta_{1}0\eta_{2}0}\left(0,;0,x_{2}\right)\nonumber \\
 & = & -\frac{\eta_{1}^{3/2}}{\sqrt{\pi}}\frac{\partial}{\partial\eta_{1}}\frac{4\pi\left(e^{-\eta_{2}x_{2}}-e^{-\eta_{1}x_{2}}\right)}{x_{2}\left(\eta_{1}^{2}-\eta_{2}^{2}\right)}\nonumber \\
 & = & \frac{8\sqrt{\pi}\eta_{1}^{5/2}}{\left(\eta_{1}^{2}-\eta_{2}^{2}\right)^{2}}\left(\frac{1}{x_{2}}e^{-\eta_{2}x_{2}}-\left(\frac{\eta_{1}^{2}-\eta_{2}^{2}}{2\eta_{1}}+\frac{1}{x_{2}}\right)e^{-\eta_{1}x_{2}}\right)\label{eq:HY}
\end{eqnarray}

One may take the limit $\eta_{2}\rightarrow\eta_{1}$ to obtain Eq.
(49) from a previous paper \cite{Stra89a} that used a very different
integration method,

\begin{eqnarray*}
S_{1s,1}^{\eta_{2}\eta_{2}0}\left(0,;0,x_{2}\right) & = & \frac{\sqrt{\pi}\left(1+x_{2}\eta_{2}\right)}{\sqrt{\eta_{2}}}e^{-\eta_{2}x_{2}}
\end{eqnarray*}

\section{Set of integral reductions for inverse powers inverse inverse binomials
in the exponentials }

\subsection{The transformation }

\noindent We may use arbitrary values of a and b in integrals of exponentials
that contain both inverse powers and $-j/(x+y)$ . We obtain

\vspace{0.9cm}

\begin{eqnarray}
R{}_{2} & \left(\right. & n\left.,\,m,\,\nu,\,a,\,b,\,c,\,h,\,j,\,0,\,0\right)=\int_{0}^{\infty}\int_{0}^{\infty}\frac{1}{y^{m/2}x^{n/2}\left(x+y\right)^{\nu/2}}f\left(\frac{xy}{x+y}\right)\nonumber \\
 & \times & e^{-\frac{a}{x}-\frac{b}{y}-c\,xy/(x+y)-h\,y/(x+y)-j/(x+y)}dx\,dy\nonumber \\
 & = & 2\int_{0}^{\infty}f(t)\,t^{-m/2}e^{-ct}\,dt\int_{0}^{\infty}s^{m+\nu-3}\left(s^{2}+t\right)^{-\frac{m}{2}-\frac{n}{2}-\nu+2}\nonumber \\
 & \times & \exp\left(-\frac{at+bs^{2}}{t\left(s^{2}+t\right)}-\frac{ht}{s^{2}+t}-\frac{js^{2}}{\left(s^{2}+t\right)^{2}}\right)ds\quad.\label{eq:R2stp0q0hj-1}
\end{eqnarray}

If we again change variables to

\begin{equation}
r=\frac{1}{s^{2}+t},\:s=\frac{\sqrt{1-rt}}{\sqrt{r}}\qquad,\label{eq:stow-1-1}
\end{equation}
 with

\begin{eqnarray}
ds & = & -\frac{2s}{\left(s^{2}+t\right)^{2}}dr\label{eq:ds2dv-1-1-1-1}
\end{eqnarray}
 we have

\begin{eqnarray}
R{}_{2} & \left(\right. & n\left.,\,m,\,\nu,\,a,\,b,\,c,\,h,\,j,\,0,\,0\right)\nonumber \\
 & = & \int_{0}^{\infty}f(t)\,t^{-m/2}e^{-\frac{b}{t}-ct}\,dt\int_{0}^{1/t}r^{\frac{n}{2}+\frac{\nu}{2}-2}(1-rt)^{\frac{m}{2}+\frac{\nu}{2}-2}\nonumber \\
 & \times & exp\left[+jr^{2}t-r(a-b+ht+j)]\right]dr\quad.\label{eq:r_int-1}
\end{eqnarray}

This is integrable for even values of $n+\nu\geq4$ if we set the
second factor in the \emph{r} integral to unity with $m=-\nu+4$ .
One may also do a binomial expansion of $(1-rt)$ for even values
of $m+\nu>4$. This is done by completing the square in the exponential
and setting 
\begin{equation}
\text{r'}=r-\frac{a-b+ht+j}{2jt}\label{eq:rpfromr}
\end{equation}
 with unit Jacobian.

\subsection{Application}

Fourier transforms of products of Yukawa or Coulomb potentials and
hydrogenic orbitals are more difficult than (\ref{eq:SYYx2}). Consider,
for instance, the case of the Fourier Transform of a product of two
Yukawa potentials centered on different positions 
\begin{equation}
S_{1}^{\eta_{1}0\eta_{2}0}\left(\mathbf{k};0,\mathbf{x}_{2}\right)=\int d^{3}x_{1}\frac{e^{-\eta_{1}x_{1}}}{x_{1}}\frac{e^{-\eta_{2}\left|\mathbf{x}_{1}-\mathbf{x}_{2}\right|}}{\left|\mathbf{x}_{1}-\mathbf{x}_{2}\right|}e^{-i\mathbf{k}\cdot\mathbf{x}_{1}}\label{eq:SYYx2p-1}
\end{equation}

The procedure is as in Section 5.2 above, but with two extra terms
in the exponential after completing the square:

\begin{eqnarray}
S_{1}^{\eta_{1}0\eta_{2}0} & \left(\right. & \mathbf{k}\left.;0,\mathbf{x}_{2}\right)=\int d^{3}x_{1}{\displaystyle \frac{1}{\sqrt{\pi}}}\nonumber \\
 & \times & \int_{0}^{\infty}\,d\rho_{1}{\displaystyle \frac{e^{-x_{1}^{2}\rho_{1}}e^{-\eta_{1}^{2}/4/\rho_{1}}}{\rho_{1}^{\;1/2}}}{\displaystyle \frac{1}{\sqrt{\pi}}}\int_{0}^{\infty}\,d\rho_{2}{\displaystyle \frac{e^{-x_{12}^{2}\rho_{2}}e^{-\eta_{2}^{2}/4/\rho_{2}}}{\rho_{2}^{\;1/2}}e^{-i\mathbf{k}\cdot\mathbf{x}_{1}}}\nonumber \\
 & = & {\displaystyle \frac{1}{\pi}}\int d^{3}x'_{1}\int_{0}^{\infty}\,d\rho_{1}{\displaystyle \frac{e^{-\eta_{1}^{2}/4/\rho_{1}}}{\rho_{1}^{\;1/2}}}\int_{0}^{\infty}\,d\rho_{2}{\displaystyle \frac{e^{-\eta_{2}^{2}/4/\rho_{2}}}{\rho_{2}^{\;1/2}}}\nonumber \\
 & \times & exp\left(-\left(\rho_{1}+\rho_{2}\right)x'{}_{1}^{2}-\frac{x_{2}^{2}\rho_{1}\rho_{2}}{\rho_{1}+\rho_{2}}-\frac{\rho_{2}i\text{ k}\cdot\mathbf{x}_{2}}{\rho_{1}+\rho_{2}}-\frac{\text{ k}^{2}}{4\left(\rho_{1}+\rho_{2}\right)}\right)\nonumber \\
 & = & \pi^{1/2}\int_{0}^{\infty}\,d\rho_{1}{\displaystyle \frac{e^{-\eta_{1}^{2}/4/\rho_{1}}}{\rho_{1}^{\;1/2}}}\int_{0}^{\infty}\,d\rho_{2}{\displaystyle \frac{e^{-\eta_{2}^{2}/4/\rho_{2}}}{\rho_{2}^{\;1/2}}}\nonumber \\
 & \times & \frac{1}{\left(\rho_{1}+\rho_{2}\right)^{3/2}}\exp\left(-\frac{x_{2}^{2}\rho_{1}\rho_{2}}{\rho_{1}+\rho_{2}}-\frac{\rho_{2}i\text{ k}\cdot\mathbf{x}_{2}}{\rho_{1}+\rho_{2}}-\frac{\text{ k}^{2}}{4\left(\rho_{1}+\rho_{2}\right)}\right)\nonumber \\
 & = & \pi^{1/2}R{}_{2}\left(4,\,4,\,0,\,\eta_{1}^{2}/4,\,\eta_{2}^{2}/4,\,x_{2}^{2},\,i\text{ k}\cdot\mathbf{x}_{2},\,\text{ k}^{2}/4,\,0,\,0\right)\label{eq:SYYx2pGaussian-2}
\end{eqnarray}
 where in the last line we have taken our function to be 
\begin{equation}
f\left(t\right)=t^{\mu}\label{eq:fispowers-1-1}
\end{equation}
 and in the present case have $\mu=3/2$. We could also have set $\mu=0$
with $n=m=1$ and $\nu=3$.

Then

\begin{eqnarray}
S_{1}^{\eta_{1}0\eta_{1}0}\  & \left(\right. & \mathbf{k}\left.;0,\mathbf{x}_{2}\right)=\pi^{1/2}R{}_{2}\left(4,\,4,\,0,\,\eta_{1}^{2}/4,\,\eta_{2}^{2}/4,\,x_{2}^{2},\,\text{i k}\cdot\mathbf{x}_{2},\,\text{k}^{2},\,0,\,0\right)\nonumber \\
 & = & \pi^{1/2}\int\,t^{-1/2}e^{-\frac{\eta_{2}^{2}}{4t}-x_{2}^{2}t}\,dt\nonumber \\
 & \times & \int_{0}^{1/t}exp\left[\text{ k}^{2}r^{2}t+r\left(-itk\cdot x_{2}-\frac{k^{2}}{4}-\frac{\eta_{1}^{2}}{4}+\frac{\eta_{2}^{2}}{4}\right)\right]dr\nonumber \\
 & = & \pi^{1/2}\int_{0}^{\infty}dt\,e^{-\frac{\eta_{2}^{2}}{4t}-x_{2}^{2}t}\frac{\sqrt{\pi}}{kt}\exp\left(-\frac{\left(i\text{ k}\cdot\mathbf{x}_{2}t+\frac{k^{2}}{4}+\frac{\eta_{1}^{2}}{4}-\frac{\eta_{2}^{2}}{4}\right)^{2}}{k^{2}t}\right)\label{eq:SYYx2pGaussian-1-1}\\
 & \times & \left(\text{erfi}\left(\frac{i\text{ k}\cdot\mathbf{x}_{2}t+\frac{k^{2}}{4}+\frac{\eta_{1}^{2}}{4}-\frac{\eta_{2}^{2}}{4}}{k\sqrt{t}}\right)-\text{erfi}\left(\frac{i\text{ k}\cdot\mathbf{x}_{2}t-\frac{k^{2}}{4}+\frac{\eta_{1}^{2}}{4}-\frac{\eta_{2}^{2}}{4}}{k\sqrt{t}}\right)\right)\nonumber 
\end{eqnarray}

As a check one can instead change variables to $\tau={\displaystyle \frac{\rho_{2}}{\rho_{1}+\rho_{2}}=\frac{y}{x+y}}$
in (\ref{eq:SYYx2pGaussian-2}) to give \cite{GR5 p. 384 No. 3.471.9}

\begin{eqnarray}
S_{1}^{\eta_{1}0\eta_{2}0} & \left(\right. & \mathbf{k}\left.;0,\mathbf{x}_{2}\right)=\pi^{1/2}\int_{0}^{1}\,d\tau{\displaystyle \frac{1}{\tau{}^{1/2}}}e^{-i\mathbf{k}\cdot\mathbf{x}_{2}\tau}\int_{0}^{\infty}\,dx{\displaystyle \frac{1}{x^{3/2}}}\label{eq:cheshire-1}\\
 & \times & \exp\left(-x_{2}^{2}\tau x-\frac{(1-\tau)\left(k^{2}\tau+\eta_{2}^{2}\right)+\eta_{1}^{2}\tau}{4\tau x}\right)\nonumber \\
 & = & \pi^{1/2}\int_{0}^{1}\,d\tau e^{-i\mathbf{k}\cdot\mathbf{x}_{2}\tau}\frac{2\sqrt{\pi}\exp\left(-x_{2}\sqrt{(1-\tau)\left(k^{2}\tau+\eta_{2}^{2}\right)+\eta_{1}^{2}\tau}\right)}{\sqrt{(1-\tau)\left(k^{2}\tau+\eta_{2}^{2}\right)+\eta_{1}^{2}\tau}}\quad,\nonumber \\
 & = & 2\sqrt{\pi}\pi^{1/2}\int_{0}^{1}\,d\tau e^{-i\mathbf{k}\cdot\mathbf{x}_{2}\tau}\frac{\exp\left(-x_{2}L\right)}{L}\nonumber 
\end{eqnarray}
 where 
\begin{equation}
L=\sqrt{(1-\tau)\left(k^{2}\tau+\eta_{2}^{2}\right)+\eta_{1}^{2}\tau}\quad.\label{eq:L}
\end{equation}
 One can show numerically that these two integrals are equal. Cheshire
\cite{Cheshire} reduced the related integral (his eq. (19))

\begin{eqnarray}
I_{1} & = & \left.\frac{\eta_{1}^{3/2}}{\sqrt{\pi}}S_{1s\,1}^{\eta_{1}0\eta_{2}0}\left(\frac{1}{2}\mathbf{k_{f}};0,x_{2}\right)\right|_{\eta_{1}=1,\eta_{2}=1/2}\nonumber \\
 & = & \left.\int d^{3}x_{1}\frac{\eta_{2}^{3/2}}{\sqrt{\pi}}e^{-\eta_{2}x_{12}}\frac{\eta_{1}^{3/2}}{\sqrt{\pi}}\frac{e^{-\eta_{1}x_{1}}}{x_{1}}e^{-i\frac{1}{2}\mathbf{k_{f}}\cdot\mathbf{x}_{1}}\right|_{\eta_{1}=1,\eta_{2}=1/2}\nonumber \\
 & = & \left.\frac{\eta_{1}^{3/2}}{\sqrt{\pi}}\frac{\eta_{2}^{3/2}}{\sqrt{\pi}}\left(-\frac{\partial}{\partial\eta_{2}}\right)\int d^{3}x_{1}\frac{e^{-\eta_{2}x_{12}}}{x_{12}}\frac{e^{-\eta_{1}x_{1}}}{x_{1}}e^{-i\frac{1}{2}\mathbf{k_{f}}\cdot\mathbf{x}_{1}}\right|_{\eta_{1}=1,\eta_{2}=1/2}\nonumber \\
 & = & \left.\frac{\eta_{1}^{3/2}}{\sqrt{\pi}}\frac{\eta_{2}^{3/2}}{\sqrt{\pi}}\left(-\frac{\partial}{\partial\eta_{2}}\right)S_{1}^{\eta_{1}0\eta_{2}0}\left(\frac{1}{2}\mathbf{k_{f}},;0,x_{2}\right)\right|_{\eta_{1}=1,\eta_{2}=1/2}\quad,\label{eq:cheshire}
\end{eqnarray}
 which matches the present result after substituting the specialized
values for $\eta_{1}$ and $\eta_{2}$.

If $\eta_{1}=\eta_{2}$ we have

\begin{equation}
L=\sqrt{k^{2}(1-\tau)\tau+\eta_{1}^{2}}\quad.\label{eq:Lred}
\end{equation}
 In the limits $k\rightarrow0$ and $\eta_{2}=\eta_{1}$ we can analytically
integrate either form to give 
\begin{equation}
S_{1}^{\eta_{1}0\eta_{1}0}\left(0,;0,x_{2}\right)=\frac{2\pi e^{-x_{2}\eta_{1}}}{\eta_{1}}\label{eq:Sred}
\end{equation}

\section{Acknowledgment}

\noindent I would like to thank Professor Ray A. Mayer, of Reed College,
for creating a very different proof of (\ref{eq:fixed}) for me, and
to his colleague Professor Nicholas Wheeler for the introduction.

\eject

{document}}
\end{document}